\documentclass[12pt]{article}
\usepackage{amssymb}
\usepackage{latexsym}

\oddsidemargin=\evensidemargin
\begin{document}
\newfont{\blb}{msbm10 scaled\magstep1}
\renewcommand{\abstractname}{ }

\newtheorem{theo}{Theorem}[section]

\newtheorem{defi}[theo]{Definition}
\newtheorem{prop}[theo]{Proposition}
\newtheorem{lemm}[theo]{Lemma}
\newtheorem{coro}[theo]{Corollary}
\pagestyle{myheadings}
\date{}
\author{Gareth Tracey (G.M.Tracey@bath.ac.uk)\\ Gunnar Traustason (G.Traustason@bath.ac.uk) \thanks{This research is supported through a standard grant from EPSRC}\\ \\
Department of  Mathematical Sciences, \\
University of Bath, \\
Bath BA2 7AY,
UK}

\title{Left $3$-Engel elements in groups of exponent $60$}
\maketitle
\begin{abstract}
\mbox{}\\
Let $G$ be a group and let $x\in G$ be a left $3$-Engel element of order dividing $60$. Suppose furthermore
that $\langle x\rangle^{G}$ has no elements of order $8$, $9$ and $25$. We show that $x$ is then contained
in the locally nilpotent radical of $G$. In particular all the left $3$-Engel elements of a group of exponent $60$ are contained in the locally nilpotent radical. 

%study left $3$-Engel elements and show in particular that they are
%contained in the Hirsh-Plotkin radical of any $n$-Engel group $G$ that
%has no element of prime order less than $n$.
%
\end{abstract}

%\noindent 2000 Mathematics Subject Classification:  20F45
%\section{Introduction}
%
\section{Introduction}
\mbox{}\\
Let $G$ be a group. An element $a\in G$ is a left Engel element in $G$, if for
each $x\in G$ there exists a non-negative integer $n(x)$ such that
      $$[[[x,\underbrace {a],a],\ldots ,a]}_{n(x)}=1.$$
If $n(x)$ is bounded above by $n$ then we say that $a$ is a left $n$-Engel element
in $G$. It is straightforward to see that any element of the Hirsch-Plotkin
radical $HP(G)$ of $G$ is a left Engel element and the converse is known
to be true for some classes of groups, including solvable groups and 
finite groups (more generally groups satisfying the maximal condition on
subgroups) [3,6]. The converse is however not true in general and this is the case
even for bounded left Engel elements. In fact whereas one sees readily that 
a left $2$-Engel element is always in the Hirsch-Plotkin radical this
is still an open question for left  $3$-Engel elements. There is some 
substantial general progress by A. Abdollahi in [1] where he proves in particular that 
for any left $3$-Engel $p$-element $a$ in a group $G$ one has that $a^{p}$
is in $HP(G)$ (in fact he proves the stronger result that $a^{p}$ is in
the Baer radical), and that the subgroup generated by two left $3$-Engel
elements is nilpotent of class at most $4$. Then in [11] it is shown that the left $3$-Engel elements in groups of exponent $5$ are in $HP(G)$. In this paper we will extend this result to groups of exponent $60$. In fact we will 
prove something quite stronger. See also [2] for some results
about left $4$-Engel elements. \\ \\
It was observed by William
Burnside [4] that every element in a group of exponent $3$  is a left $2$-Engel
element and so the fact that every left $2$-Engel element lies in the Hirsch-Plotkin radical can be seen as the underlying reason why groups of exponent
$3$ are locally finite. For groups of $2$-power exponent there is a close link
with left Engel elements. If  $G$ is a group of exponent 
$2^{n}$ then it is not difficult to see that any element $a$ in $G$ of order $2$ is a left $(n+1)$-Engel element of $G$ (see the introduction of [11] for details). For sufficiently large $n$ we know that the variety of groups of exponent $2^{n}$ is not locally finite [8,9]. As a result one can see [11] that it follows that for sufficently large $n$ we do not have in general that a left $n$-Engel element is contained in the Hirsch-Plotkin radical. Using the fact that groups of exponent $4$ are locally finite [10], one can also see that if all left $4$-Engel elements of a group $G$ of exponent $8$ are in
$HP(G)$ then $G$ is locally finite. \\ \\
In this paper we continue our study of left $3$-Engel elements started in [11]. We first make the observation that an element $a\in G$ is a left $3$-Engel element if and only
if $\langle a,a^{x}\rangle$ is nilpotent of class at most $2$ for all 
$x\in G$ [1]. In [11] we introduced the following related class of groups. \\ \\
{\bf Definition}. A {\it sandwich} group is a group $G$ generated
by a  set $X$ of elements such that $\langle x,y^{g}\rangle$ is nilpotent
of class at most $2$ for all $x,y\in X$ and all $g\in G$. \\ \\
{\bf Remark}. In [11] it was shown that any sandwich group of rank $3$ is nilpotent. \\ \\
If $a\in G$ is a left $3$-Engel element then $H=\langle a\rangle^{G}$
is a {\it sandwich} group and it is clear that the following statements
are equivalent: \\ \\
(1) For every pair $(G,a)$ where $a$ is a left $3$-Engel element in the
group $G$
we have that $a$ is in the locally nilpotent radical of $G$. \\ \\
(2) Every sandwich group is locally nilpotent. \\ \\
It is also clear that to prove (2), it suffices to show that every
finitely generated sandwich group is nilpotent. \\ \\
{\bf Left $3$-Engel elements of finite order}. For left $3$-Engel elements of finite order some further reduction can be made. Suppose $G$ is a group with a left $3$-Engel element
$x$ of order $m=p_{1}^{n_{1}}\cdots p_{r}^{n_{r}}$ where $p_{1},\ldots ,p_{r}$ are distinct primes and $n_{1},\ldots ,n_{r}$ are positive integers. For $1\leq j\leq r$, let
$m_{j}=m/p_{j}^{n_{j}}$. Then $m_{1},\ldots ,m_{r}$ are coprime. Thus in order to show that $x\in HP(G)$, it suffices to show that $x^{m_{1}},\ldots ,x^{m_{r}}\in HP(G)$.
So the problem of showing that an element of finite order is in $HP(G)$ reduces to dealing with elements of prime power order. Further reductions can be made.  First we recall
a standard notion. Let $G$ be a group. For any set $\pi$ consisting of primes, we say that $x$ is a $\pi$-element in $G$ if the order of $G$ only has numbers from $\pi$ as prime
factors. 
\begin{lemm} Let  $\pi$ be a set of primes. Suppose that for all groups $G$ and all primes $p\in \pi$ we have that all left $3$-Engel elements of order $p$ in $G$ are contained
in $HP(G)$. It then follows that for all groups $G$, all left $3$-Engel elements in $G$ that are $\pi$-elements are in $HP(G)$. 
\end{lemm}
{\bf Proof}\  \ Let $G$ be any group. We have already seen that we only need to consider the case when $x$ is a left $3$-Engel element in $G$ of some prime power exponent $p^{n}$ where
$p\in \pi$ and $n$ is a positive integer. By the result of Alireza [1] mentioned above, we know that $x^{p}\in HP(G)$. As $x^{p}$ is a $p$-element we then know
that $N=\langle x^{p}\rangle^{G}$ is a locally finite $p$-group. By our assumption we know that $\langle x\rangle^{G} N/N$ is locally nilpotent and thus a locally
finite $p$-group. Hence $\langle x\rangle^{G}$ is a locally finite $p$-group and thus locally nilpotent. We thus conclude that $x\in HP(G)$. $\Box$ \\ \\
{\bf Remark}. The problem of showing that all left $3$-Engel elements of finite order are in the Hirsch-Plotkin radical thus reduces to only having to consider elements
of prime order. Thus dealing with left $3$-Engel elements of finite order reduces to working with sandwich groups generated by elements of prime order $p$. This is because the
following are equivalent for any prime $p$:  \\ \\
($1_{p}$) For every pair $(G,a)$ where $a$ is a left $3$-Engel element of order $p$ in $G$ we have that $a\in HP(G)$. \\ \\
($2_{p}$) Every finitely generated sandwich group generated by elements of order $p$ is nilpotent. \\ \\
In sections 2 and 3. We will work with sandwich groups of rank 3 and 4 generated by elements of order $2$. \\ \\
{\bf Left $3$-Engel elements in groups of finite exponent}. Determining whether a left $3$-Engel element of finite order is in the Hirsch-Plotkin radical seems a very difficult problem in general. One could thus consider adding further constraints on the group. For example one could require that for the given left $3$-Engel element $x$ in $G$ we have that
$\langle x\rangle^{G}$ is of finite exponent. In fact we will consider a weaker condition.  Let $p_{1}^{n_{1}},\ldots ,p_{r}^{n_{r}}$ be non-trivial powers where 
$p_{1},\ldots ,p_{r}$ are distinct primes. Consider the following statement. 
\begin{eqnarray*}
    {\mathcal E}(p_{1}^{n_{1}},\ldots ,p_{r}^{n_{r}}): &   &  \mbox{For all groups }G\mbox{ and all left }3\mbox{-Engel elements }x\in G\mbox{ of order } \\
     &  & \mbox{dividing }
p_{1}^{n_{1}}\cdots p_{r}^{n_{r}},\mbox{ where }\langle x\rangle^{G}\mbox{ has no elements of order }p_{1}^{n_{1}+1}, \\
      &  & \ldots ,p_{r}^{n_{r}+1}, \mbox{ we
have that }x\in HP(G).
\end{eqnarray*}
{\bf Remark}. Notice that if ${\mathcal E}(p_{1}^{n_{1}},\ldots ,p_{r}^{n_{r}})$ holds, then it would follow that a left $3$-Engel element of $G$ is in
$HP(G)$ when $\langle x\rangle^{G}$ has exponent dividing $p_{1}^{n_{1}}\cdots p_{r}^{n_{r}}$. Thus in particular all left $3$-Engel elements in a group
$G$ of exponent dividing $p_{1}^{n_{1}}\cdots p_{r}^{n_{r}}$ would be in the Hirsch-Plotkin radical. \\ \\
We will next prove a reduction result that is similar in nature to Lemma 1.1. For a given prime $p$ and positive integer $n$, consider the following statement 
\begin{eqnarray*}
     {\mathcal Q}(p,n): & & \mbox{For all groups }G\mbox{ and all  left }3\mbox{-Engel elements }x\in G\mbox{, where }x\mbox{ is of order }p \\
		& & \mbox{ and }\langle
x\rangle^{G}\mbox{ has no element of order }p^{n+1},\mbox{ we have that }x\mbox{ is in }\mbox{HP}(G). 
\end{eqnarray*}
\begin{prop} Let  $m=p_{1}^{n_{1}}\cdots p_{r}^{n_{r}}$ be an integer where $p_{1},\ldots ,p_{r}$ are distinct primes and $n_{1},\ldots ,
n_{r}$ are positive integers. Then
  $${\mathcal Q}(p_{1},n_{1})\wedge
	\cdots \wedge {\mathcal Q}(p_{r},n_{r})
	\Rightarrow {\mathcal E}(p_{1}^{n_{1}},
	\ldots ,p_{r}^{n_{r}}).$$
\end{prop}
{\bf Proof}\  \   Let $G$ be any group and $x$ a left $3$-Engel element in $G$ of order dividing $m$. Suppose that $\langle x\rangle^{G}$ has no elements of order $p_{1}^{n_{1}+1},\ldots ,p_{r}^{n_{r}+1}$. For $1\leq j\leq r$, let
$m_{j}=m/p_{j}^{n_{j}}$. Then $m_{1},\ldots ,m_{r}$ are coprime. Thus in order to show that $x$ is in $\mbox{HP}(G)$,  it suffices to show
that $x^{m_{1}},\ldots , x^{m_{r}}$ are in $\mbox{HP}(G)$. It thus suffices to deal with the case when  $m=p^{n}$ for a prime $p=p_{j}$ and
the positive integer $n=n_{j}$. Let $x_{1},\ldots ,x_{r}$ be finitely many conjugates of $x$. We want to show that $H=\langle x_{1},\ldots ,x_{r}\rangle$ is nilpotent. By the result 
of A. Abdollahi [1] mentioned above, we know that $N=\langle x_{1}^{p},\ldots ,x_{r}^{p}\rangle^{H}$ is locally nilpotent. As any finitely generated subgroup of $N$ is contained in a subgroup of $N$ generated by finitely many conjugates of $x^{p}$ and as $x^{p}$ is of $p$-power order, it follows that $N$ is a $p$-group. As $N$ is a $p$-group and $H$
contains no elements of order $p^{n+1}$, the same is true for $H/N$. As ${\mathcal Q}(p,n)$ holds by assumption, we thus have that  $\langle x_{i}\rangle^{H}N/N$ is locally
nilpotent for $1\leq i\leq r$ and hence $H/N$ is nilpotent. Thus $H/N$ is a finite $p$-group that implies that $N$ is
finitely generated and thus also a finite $p$-group. We conclude that $H$ is a finite $p$-group and thus nilpotent. We have thus shown that $\langle x\rangle^{G}$ is locally nilpotent and therefore that $x\in HP(G)$. $\Box$ \\ \\
{\bf Remark}. In particular, it follows from last proposition that in order to show that left $3$-Engel elements in groups of finite exponent are in the Hirsch-Plotkin radical, it suffices to show that ${\mathcal Q}(p,n)$ holds for all primes and positive integers $n$. \\ \\
At this stage we don't even know if for a group $G$, satisfying the hypothesis in ${\mathcal Q}(p,n)$, we can conclude that $\langle x\rangle^{G}$ is a $p$-group. Let us thus consider the following weaker statement. 
\begin{eqnarray*}
     {\mathcal R}_{G}(p,n): & & \mbox{For all  left }3\mbox{-Engel elements }x\in G\mbox{, where }x\mbox{ is of order }p\mbox{ and } \\
		& & \langle
x\rangle^{G}\mbox{ has no element of order }p^{n+1},\mbox{ we have that }\langle x\rangle^{G}\mbox{ is a }p\mbox{-group}. 
\end{eqnarray*}
{\bf Remark}\ Of course ${\mathcal R}_{G}(p,n)$ implies that $\langle x\rangle^{G}$ is of exponent dividing $p^{n}$. The next proposition gives us a sufficient condition for ${\mathcal R}_{G}(p,n)$.  
\begin{prop} Let $p$ be a prime and $n$ a positive integer. Let $G$ be a group with the property that for any left $3$-Engel element $x\in G$ of order $p$ we have that 
$\langle x\rangle^{G}$ has no element of order $p^{n+1}$. Now suppose furthermore that for any left $3$-Engel element $x$ of order $p$ and any element $g\in \langle x\rangle^{G}$ of order dividing $p^{n}$ we have
    $$\langle x,x^{g},\ldots ,x^{g^{p^{n}-1}}\rangle$$
is nilpotent. Then ${\mathcal R}_{G}(p,n)$ holds. 
\end{prop}
{\bf Proof}\ \ Let $y\in \langle x\rangle^{G}$. Then $y=x_{1}\cdots x_{r}$ for some $r$ conjugates of
$x$. We show by induction on $r$ that $y^{p^{n}}=1$. This is obvious when $r=0$. Now let $r\geq 1$ and suppose that our claim holds for smaller values of $r$. By the induction hypothesis
we know that $(x_{1}\cdots x_{r-1})^{p^{n}}=1$. Thus for $m=p^{n}$ and $g=x_{1}\cdots x_{r-1}$ we have
\begin{eqnarray*}
   y^{m} & = & g^{m}x_{r}^{g^{m-1}}x_{r}^{g^{m-2}}\cdots x_{r}^{g}x_{r} \\
          & = & x_{r}^{g^{m-1}}\cdots x_{r}^{g}x_{r}.
\end{eqnarray*}
By our assumptions $\langle x_{r},x_{r}^{g},\ldots ,x_{r}^{g^{m-1}}\rangle$ is nilpotent and thus a finite
$p$-group. By the assumptions $\langle x_{r},
x_{r}^{g},\ldots ,x_{r}^{g^{m-1}}\rangle$ is then of exponent $m$. In particular it follows that 
$$(x_{r}^{g^{m-1}}\cdots x_{r}^{g}x_{r})^{m}=1$$ 
and thus $y^{m^{2}}=1$. As $\langle x\rangle^{G}$ has no
element of order $pm$ it follows that $y^{m}=1$. $\Box$ \\ \\
{\bf Remark}. As any sandwich group of rank $3$ is nilpotent [11], it follows from Proposition 1.3 that ${\mathcal R}_{G}(2,1)$ and ${\mathcal R}_{G}(3,1)$ hold in any group $G$. Groups of exponent $2$ are abelian and from Burnside [4] we know that groups of exponent $3$ are locally finite. It thus follows that ${\mathcal Q}(2,1)$ and ${\mathcal Q}(3,1)$ hold. \\ \\
The main result of the paper is the following. 
\begin{theo} Let $G$ be any group and let $x$ be a left $3$-Engel in $G$ of order dividing $60$. Suppose furthermore
that $\langle x\rangle^{G}$ has no elements of order $8,9$ or $25$. Then $x\in HP(G)$. 
\end{theo}
{\bf Remark}. From Proposition 1.2 it suffices to show that ${\mathcal Q}(2,2), {\mathcal Q}(3,1)$
and ${\mathcal Q}(5,1)$ hold. We have already seen from last remark that ${\mathcal Q}(3,1)$ holds. It thus
remains to see that ${\mathcal Q}(2,2)$ and ${\mathcal Q}(5,1)$ hold. \\ \\
To prove ${\mathcal Q}(2,2)$ we will need some preliminary work. This will be carried out in Sections 2 and 3. In [11] it was shown that all sandwich groups of rank $3$ are nilpotent. The proof for the case when the
group is generated by involutions is substantially simpler and thus we start by giving a short proof of this in Section 2. In Section 3 we will then deal with certain sandwich groups of rank $4$ generated by involutions that are needed to prove ${\mathcal Q}(2,2)$. In Sections 4 and 5 we will the prove ${\mathcal Q}(2,2)$ and ${\mathcal Q}(5,1)$ respectively. \\ \\
{\bf Remark}. Our way of writing polycyclic presentations in this paper follows [7]. It reflects a polycyclic series 
   $$\langle x_{1}\rangle\unlhd
	\langle x_{1},x_{2}\rangle\unlhd
	\cdots \unlhd \langle x_{1},\ldots ,x_{m}\rangle =G.$$
We also partition the set of generators into subsets $X_{1},\ldots ,X_{r}$ where
$\langle X_{1}\rangle\leq \langle X_{1}\cup
X_{2}\rangle\leq \cdots \langle X_{1}\cup
\cdots \cup X_{r}\rangle =G$ is a normal series with abelian factors. 
\section{Sandwich groups generated by $3$ involutions}
Let $F=\langle x,y,z\rangle$ be a $3$-generator sandwich group generated by involutions $x,y$ and $z$. 
\begin{theo} $F$ is nilpotent of class at most $5$. 
\end{theo}
{\bf Proof}\ \ We have that $a=x$ and $b=x^{y}$ commute. Thus $a,b^{z}$ commute with $a^{z},b$ and 
    $$\langle a,b,z\rangle = (\langle a,b^{z}\rangle\cdot \langle a^{z},b\rangle)\ltimes \langle z\rangle.$$
Then 
\begin{eqnarray*}
  \mbox{}[x,y,z,z] & = & [ab,z,z] \\
                              & = & [(ab^{z})(ba^{z}),z] \\
                             & = & (b^{z}a)^{2}(a^{z}b)^{2} \\
                             & = & [b^{z},a]\cdot [a^{z},b].
\end{eqnarray*}
This element clearly commutes with $a=x$ and $z$. As $\langle x,y\rangle$ is nilpotent of class at most $2$, we have $1=[x^{2},y]=[x,y]^{2}$ and
thus $[x,y,z,z]=[y,x,z,z]$. By symmetry we thus see that $[x,y,z,z]=[y,x,z,z]$ commutes with $y$ and is thus in $Z(F)$. By symmetry it follows that 
\begin{equation}
                           [x,y,z,z]=[y,x,z,z],\,[y,z,x,x]=[z,y,x,x],\,[z,x,y,y]=[x,z,y,y]\in Z(F).
\end{equation}
Next notice that 
     $$[[x,y],[x,z]]=[ab,aa^{z}]=[ab,a^{z}]=[b,a^{z}]$$
commutes with $a=x$ and we have seen in (1) that modulo $Z(F)$ we have $[b,a^{z}]=[b^{z},a]$. Thus, modulo $Z(F)$, we know that
$[[x,y],[x,z]]$ commutes with $x$ and $z$. By symmetry we then see that  $[[x,y],[x,z]]=[[x,z],[x,y]]^{-1}$ also commutes with $y$ modulo $Z(F)$. Hence
$[[x,y],[x,z]]\in Z_{2}(F)$. By symmetry
\begin{equation}
          [[x,y],[x,z]],\,[[y,z],[y,x]],\,[[z,x],[z,y]]\in Z_{2}(F).
\end{equation}
Then, modulo $Z_{2}(F)$, we have 
\begin{eqnarray*}
          [x,y,z]^{x} & = & [x,y,z[z,x]] \\
                & = & [x,y,[z,x]z] \\
                & = & [x,y,z]\cdot [x,y,[z,x]]^{z} \\
               & = & [x,y,z].
\end{eqnarray*}
Thus from (1) and (2) we know that $[x,y,z]=[y,x,z]$ commutes with $x,y,z$ modulo $Z_{2}(F)$. By symmetry we thus have
\begin{equation}
[x,y,z],[y,z,x],[z,x,y]\in Z_{3}(F).
\end{equation}
As $\langle x,y\rangle,\langle y,z\rangle ,\langle z,x\rangle$ are nilpotent of class at most $2$, it follows from (3) that $[x,y],[y,z],[z,x]\in Z_{4}(F)$ from
which it follows that $x,y,z\in Z_{5}(F)$. $\Box$ \\ \\
Knowing that $F$ is nilpotent, it is now easy to come up with a power-conjugation presentation for the largest such group. Notice first that
        $$[z,x,y,[z,x]]=[z,x,y,z,x][z,x,y,x,z]=[z,x,[y,z],x][x,z,[y,x],z].$$
Calculating in $\langle a,b,x\rangle=\langle z,z^{y},x\rangle $, we see that $[z,x,[y,z],x]=[a,x,ab,x]=[a^{x},b,x]=[a^{x},b][a,b^{x}]= [z,y,x,x]$. By symmetry
we have $[x,z,[y,x],z]=[x,y,z,z]$. Then, calculating in $\langle a,b,y\rangle =\langle z,z^{x},y\rangle$, we see that $[z,x,y,[z,x]]=[ab,y,ab]=
[a^{y}b^{y},ab]=[a^{y},b][b^{y},a]=[z,x,y,y]$. We thus have $[z,x,y,y]=[z,y,x,x][x,y,z,z]$ or
                            $$[z,x,y,y][x,y,z,z][y,z,x,x]=1.$$
One can come up with a full presentation using for example the nilpotent quotient algorithm or by hand. It turns our that we get the group $F=\langle
x,y,z\rangle$ of order $2^{13}$ where the generators and relations are as follows. \\ \\
\underline{Generators} \\ \\
$X_{1}:\ \ x_{1}=[z,x,y,y],\ x_{2}=[x,y,z,z],\
x_{3}=[y,z,x,x]$ \\ \\
$X_{2}:\ \ x_{4}=[z,x,[z,y]],\ x_{5}=[x,y,[x,z]],\
x_{6}=[y,z,[y,x]],$ \\
$\mbox{\ \ \ \ \ \ \ \ }x_{7}=[z,x,y],\ x_{8}=[z,y,x]$ \\ \\
$X_{3}:\ \ x_{9}=[z,x],\ x_{10}=[z,y],\ x_{11}=[x,y]$ \\ \\
$X_{4}:\ \ x_{12}=x,\  x_{13}=y,\ x_{14}=z.$ \\ \\
\underline{Relations} \\ \\
$x_{3}=x_{2}x_{1}$, \\ 
$x_{1}^{2}=x_{2}^{2}=x_{4}^{2}=x_{5}^{2}=x_{6}^{2}=x_{9}^{2}=x_{10}^{2}=x_{11}^{2}=x_{12}^{2}=x_{13}^{2}=x_{14}^{2}=1,\ 
x_{7}^{2}=x_{1}, x_{8}^{2}=x_{2}x_{1} $,\\ \\
$x_{4}^{x_{12}}=x_{4}x_{2}x_{1},\ x_{4}^{x_{13}}=x_{4}x_{1},\
x_{5}^{x_{13}}=x_{5}x_{1},\ x_{5}^{x_{14}}=x_{5}x_{2}$,\\
$x_{6}^{x_{12}}=x_{6}x_{2}x_{1},\ x_{6}^{x_{14}}=x_{6}x_{2},$ \\
$x_{7}^{x_{9}}=x_{7}x_{1},\ x_{7}^{x_{10}}=x_{7}x_{1},\
x_{7}^{x_{11}}=x_{7}x_{1},\ x_{7}^{x_{12}}=x_{7}x_{5}x_{1},$ \\
$x_{7}^{x_{13}}=x_{7}x_{1},\ x_{7}^{x_{14}}=x_{7}x_{4}x_{1},$
\\
$x_{8}^{x_{9}}=x_{8}x_{2}x_{1},\ x_{8}^{x_{10}}=x_{8}x_{2}x_{1},\
x_{8}^{x_{11}}=x_{8}x_{2}x_{1},\ x_{8}^{x_{12}}=x_{8}x_{2}x_{1},$ \\
$x_{8}^{x_{13}}=x_{8}x_{6}x_{2}x_{1},\ x_{8}^{x_{14}}=x_{8}x_{4}x_{2}x_{1},$
\\
$x_{9}^{x_{10}}=x_{9}x_{4},\ x_{9}^{x_{11}}=x_{9}x_{5},\
x_{9}^{x_{13}}=x_{9}x_{7},\ x_{10}^{x_{11}}=x_{10}x_{6},\
x_{10}^{x_{12}}=x_{10}x_{8},$ \\
$x_{11}^{x_{14}}=x_{11}x_{8}x_{7}x_{6}x_{5}x_{4}x_{2},\
x_{12}^{x_{13}}=x_{12}x_{11},\ x_{12}^{x_{14}}=x_{12}x_{9},\
x_{13}^{x_{14}}=x_{13}x_{10}$.
\section{Sandwich groups generated by $4$ involutions}
In this section we move on to $4$-generator sandwiches. The ultimate aim is to show
that these are nilpotent. We get here some partial results that will be sufficient to prove the main results of this paper. This is achieved by analysing various quotients
of the largest sandwich group of rank $4$ generated by involutions. The following definition will be useful. \\ \\
{\bf Definition}. Let $G$ be a sandwich group generated by a finite set
$X=\langle a_{1},\ldots ,a_{r}\rangle$ of sandwich elements. The {\it commutativity graph} of $G$, $V(G)$, is an (undirected) graph whose set of vertices 
is the set of generators $X$ and where a pair of distinct vertices $a_{i}$ and 
$a_{j}$ are joined by and edge if and only if $a_{i}$ and $a_{j}$ commute. \\ \\
{\bf Remarks}. (1) The commutativity graph of the free $r$-generator sandwich
group has no edges and the largest $r$-generator sandwich group, whose 
commutativity graph is the complete graph, is the free abelian group of
rank $r$. \\ \\
(2) Let $H$ and $K$ be the largest $r$-generator sandwich groups with 
commutativity graphs $V(H)$ and $V(K)$ respectively. If $V(H)\subseteq V(K)$
then $K$ is isomorphic to a quotient of $H$. \\ \\
We now focus on sandwich groups generated by $4$ involutions. It is clear
that if we have a complete commutativity graph we get  $C_{2}^{4}$ that is of order $16$. There is only one type of
a commutativity graph with $5$ edges, namely

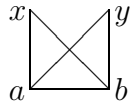
\mbox{}\\ \\
\begin{picture}(100,50)(-120,-50)
%
%
%\put(15,5){\line(1,0){30}}
\put(15,-25){\line(1,0){30}}
%\put(30,20){\line(3,5){30}}
%\put(90,20){\line(-3,5){30}}
\put(47,-29){$b$}
\put(7,-29){$a$}
\put(47,1){$y$}
\put(7,1){$x$}
\put(15,5){\line(1,-1){30}}
\put(15,-25){\line(1,1){30}}
\put(15,5){\line(0,-1){30}}
\put(45,5){\line(0,-1){30}}
\end{picture}
\normalsize
\mbox{}\\
and the largest $4$-generator sandwich group $\langle x,y,a,b\rangle$
with this commutativity graph, that is generated by involutions, is $\langle x,y\rangle\times \langle a,b\rangle=
D_{8}\times C_{2}^{2}$ that is of order $32$.  Moving next on to sandwich groups, whose commutativity graph has
$4$ edges, there are the following two types of graphs to conisder (either
the two removed edges are adjacent or not) \\
 \mbox{}\\ \\
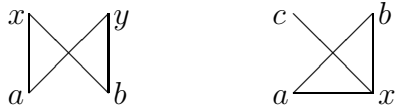
\begin{picture}(100,50)(-120,-50)
%
%
%\put(-5,5){\line(1,0){30}}
%\put(-5,-25){\line(1,0){30}}
%\put(30,20){\line(3,5){30}}
%\put(90,20){\line(-3,5){30}}
\put(27,-29){$b$}
\put(-13,-29){$a$}
\put(27,1){$y$}
\put(-13,1){$x$}
\put(-5,5){\line(1,-1){30}}
\put(-5,-25){\line(1,1){30}}
\put(-5,5){\line(0,-1){30}}
\put(25,5){\line(0,-1){30}}
%
%\put(135,-25){\line(1,0){30}}
\put(95,-25){\line(1,0){30}}
%\put(30,20){\line(3,5){30}}
%\put(90,20){\line(-3,5){30}}
\put(127,-29){$x$}
\put(87,-29){$a$}
\put(127,1){$b$}
\put(87,1){$c$}
\put(95,5){\line(1,-1){30}}
\put(95,-25){\line(1,1){30}}
%\put(95,5){\line(0,-1){30}}
\put(125,5){\line(0,-1){30}}

\end{picture}
\normalsize
\mbox{}\\
The largest sandwich group where $x,y,a,b$ are involutions and with the former commutivity graph is $\langle x,y\rangle\times
\langle a,b\rangle=D_{8}\times D_{8}$ that has order $64$. Moving to the latter group notice that $a,b^{c}$ commute with $b,a^{c}$ and we thus have
     $$\langle a,b,c\rangle=\langle a,b^{c}\rangle \wr \langle c\rangle = D_{8}\wr C_{2}$$ 
that is the standard wreath product of $D_{8}$ by $C_{2}$. Thus the largest sandwich group generated by involutions $a,b,c,x$ that have the latter
graph as a commutivity graph is 
    $$\langle a,b,c\rangle\times \langle x\rangle=(D_{8}\wr C_{2})\times C_{2}$$
and is of order $256$. \\ \\
Next we consider the case when the commutativity graph has 3 edges. This is 
much more difficult and needs some care. Here there are three types
of commutativity graphs. These are \\
 \mbox{}\\ \\
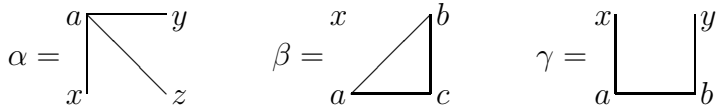
\begin{picture}(100,50)(-120,-50)
\put(-55,5){\line(1,0){30}}
%\put(-55,-25){\line(1,0){30}}
%\put(30,20){\line(3,5){30}}
%\put(90,20){\line(-3,5){30}}
\put(-23,-29){$z$}
\put(-63,-29){$x$}
\put(-23,1){$y$}
\put(-63,1){$a$}
\put(-85,-14){$\alpha=$}
\put(-55,5){\line(1,-1){30}}
%\put(-55,-25){\line(1,1){30}}
\put(-55,5){\line(0,-1){30}}
%\put(-25,5){\line(0,-1){30}}
%
%\put(45,5){\line(1,0){30}}
\put(45,-25){\line(1,0){30}}
%\put(30,20){\line(3,5){30}}
%\put(90,20){\line(-3,5){30}}
\put(77,-29){$c$}
\put(37,-29){$a$}
\put(15,-14){$\beta=$}
\put(77,1){$b$}
\put(37,1){$x$}
%\put(45,5){\line(1,-1){30}}
\put(45,-25){\line(1,1){30}}
%\put(45,5){\line(0,-1){30}}
\put(75,5){\line(0,-1){30}}
%
%\put(145,5){\line(1,0){30}}
\put(145,-25){\line(1,0){30}}
%\put(30,20){\line(3,5){30}}
%\put(90,20){\line(-3,5){30}}
\put(177,-29){$b$}
\put(137,-29){$a$}
\put(115,-14){$\gamma=$}
\put(177,1){$y$}
\put(137,1){$x$}
%\put(145,5){\line(1,-1){30}}
%\put(145,-25){\line(1,1){30}}
\put(145,5){\line(0,-1){30}}
\put(175,5){\line(0,-1){30}}
\end{picture}
\normalsize
\mbox{}\\
The largest sandwich group with the first commutativity graph is
    $$G_{\alpha}=\langle x,y,z\rangle\times \langle a\rangle =R\times C$$
where $R$ is the largest sandwich group generated by involutions $x,y,z$ that we dealt with
in last section. In the next two subsections we deal with the other two types.
\subsection{Sandwich groups with commutativity graph $\beta$}
Let $G_{\beta}=\langle x,a,b,c\rangle$ be a sandwich group where $x,a,b,c$ are involutions and whose commutativity graph
is
 \mbox{}\\ \\ \\
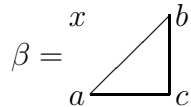
\begin{picture}(100,50)(-120,-50)
 \put(45,-25){\line(1,0){30}}
%\put(30,20){\line(3,5){30}}
%\put(90,20){\line(-3,5){30}}
\put(77,-29){$c$}
\put(37,-29){$a$}
\put(15,-14){$\beta=$}
\put(77,1){$b$}
\put(37,1){$x$}
%\put(45,5){\line(1,-1){30}}
\put(45,-25){\line(1,1){30}}
%\put(45,5){\line(0,-1){30}}
\put(75,5){\line(0,-1){30}}
\end{picture}
\normalsize
\mbox{}\\
In this subsection we show that this group is nilpotent and obtain a consistent
presentation for the group. We will use the fact that $3$-generator sandwich groups are nilpotent. The following subgroups generated by $3$ involutions will play a key role in
the following: $H(c)=\langle x,x^{ab},c\rangle,\,H(a)=\langle
x,x^{bc},a\rangle,\,H(b)=\langle x,x^{ca},b\rangle,\,K=\langle
x^{a},x^{b},x^{c}\rangle$. 
\begin{lemm} We have that \\ \\
(1)\ $[x,x^{ab}]=[x^{a},x^{b}]$ commutes with $x,a,b$.\\
(2)\ $[x,x^{bc}]=[x^{b},x^{c}]$ commutes with $x,b,c$. \\
(3)\ $[x,x^{ca}]=[x^{c},x^{a}]$ commutes with $x,c,a$. 
\end{lemm}
{\bf Proof}\ \ We have that 
\begin{eqnarray*}
\mbox{}[x,x^{ab}] & = & [x,[x,ab]] \\
                   & = & [x,(a^{x}b)(b^{x}a)] \\
               & = & (b^{x}a)(a^{x}b)(a^{x}b)(b^{x}a) \\
                & = & [b^{x},a][a^{x},b]
\end{eqnarray*}
and as $\langle b^{x},a\rangle$, $\langle a^{x},b\rangle$ are nilpotent of class at most $2$, it follows that $[x,x^{ab}]$ commutes with $a,b$. As $[x,x^{ab}]$ clearly commutes with $x$ we have that the first part follows and thus the others by symmetry. $\Box$ 
\begin{lemm} $\gamma_{5}(K)=\{1\}$.
\end{lemm}
{\bf Proof}\ \ Using the fact that $x$ commutes with $x^{a},x^{b}$ and $x^{c}$ and that by Theorem 2.1 we know that  
$H(c)$ is nilpotent of class at most $5$, we have				
\begin{eqnarray*}
\mbox{}[x^{a},x^{b},x^{c},x^{c}] & = & [x^{a},x^{b},[x,c],x^{c}] \\
 					& = & [x^{a},x^{b},[x,c],[x,c]] \\
          & \stackrel{\tiny \mbox{L3.1}}{=} & [x,x^{ab},[x,c],[x,c]] \\
					 & = & 1.
\end{eqnarray*}
By symmetry $[x^{b},x^{c},x^{a},x^{a}]=[x^{c},x^{a},x^{b},x^{b}]=1$ and from the work on the $3$-generator sandwich groups we know that implies that $\gamma_{5}(K)=\{1\}$. $\Box$														
\begin{lemm} $\gamma_{5}(H(a)),\,\gamma_{5}(H(b)),\,\gamma_{5}(H(c))\leq Z(\langle x\rangle^{G_{\beta}})$.
\end{lemm}								
{\bf Proof}\ \ From our analysis of the $3$-generator groups we know that $\gamma_{5}(H(c))$ is generated
by $[x,c,x^{ab},x^{ab}]$ and $[x,x^{ab},c,c]$ and that $[x^{ab},c,x,x]=[x,c,x^{ab},x^{ab}][x,x^{ab},c,c]$. As $H(c)$ is nilpotent of class at most $5$, we also know that these commutators all commute with $x,x^{ab}$ and $c$. By Lemma 3.1 we also know that $[x,x^{ab}]$ commutes with $a,b$ and thus $[x,x^{ab},c,c]$ commutes with $x,a,b,c$ and is thus
in $Z(G_{\beta})$. Hence 
                $$\gamma_{5}(H(c))Z(G_{\beta})=\langle [x,c,x^{ab},x^{ab}]\rangle Z(G_{\beta}).$$
In order to finish the proof of the lemma, it thus suffices to show that $[x,c,x^{ab},x^{ab}]$ commutes
with $x,x^{a},x^{b},x^{c},x^{ab},x^{ac},x^{bc},x^{abc}$. As we know already that it commutes with $x,x^{ab},c$ it 
only remains to show that it commutes with $x^{a},x^{b}$. Now using the fact again that $H(c)$ is nilpotent
of class at most $5$ and that $[x,x^{ab},x^{ab}]=[x^{c},x^{ab},x^{ab}]=1$ we have
\begin{eqnarray*}
 [x,c,x^{ab},x^{ab}] & = & [xx^{c},x^{ab},x^{ab}] \\
                & = & [[x,x^{ab}][x,x^{ab},x^{c}][x^{c},x^{ab}],x^{ab}] \\
         & = & [x,x^{ab},x^{c},x^{ab}].
\end{eqnarray*}
As $ab=ba$ we thus see from the symmetry that we only need to
show that $[x,x^{ab},x^{c},x^{ab}]$ commutes with $x^{a}$. This
will follow from the following calculations. We use there the
fact from Lemma 3.1 that $[x,x^{ab}]=[x^{a},x^{b}]$ and also
that $K=\langle x^{a},x^{b},x^{c}\rangle$ is nilpotent of class
at most $4$. We have
\begin{eqnarray*}
  [x,x^{ab},x^{c},x^{ab}]^{x^{a}} & = & [[x^{a},x^{b},x^{c}]^{x^{a}},x^{ab}] \\
        & = & [x^{a},x^{b},x^{c}[x^{c},x^{a}],x^{ab}] \\
        & = & [[x^{a},x^{b},[x^{c},x^{a}]][x^{a},x^{b},x^{c}],x^{ab}].
\end{eqnarray*}
Now we kow by Lemma 3.1 that $[[x^{a},x^{b}],[x^{c},x^{a}]]$ commutes with $a$ and as $K$ is nilpotent
of class at most $4$, it follows that 
      $$1=[x^{a},x^{b},[x^{c},x^{a}],x^{b}]^{a}=[x^{a},x^{b},[x^{c},x^{a}],x^{ab}].$$
From this and the calculations above it thus follows that 
   $$[x,x^{ab},x^{c},x^{ab}]^{x^{a}}=[[x^{a},x^{b},[x^{c},x^{a}]][x^{a},x^{b},x^{c}],x^{ab}]=[x^{a},x^{b},x^{c},x^{ab}]=[x,x^{ab},x^{c},x^{ab}],$$
Thus $[x,x^{ab},x^{c},x^{ab}]$ commutes with $x^{a}$ and this finishes the proof. $\Box$ \\
\begin{lemm} We have that the following identities hold modulo $Z(\langle x\rangle^{G_{\beta}})$. \\ \\
(1) $[x,x^{abc},x^{a}]=[x^{b},x^{c},x^{a}]$.\\
(2) $[x,x^{abc},x^{b}]=[x^{c},x^{a},x^{b}]$. \\
(3) $[x,x^{abc},x^{c}]=[x^{a},x^{b},x^{c}]$.
\end{lemm}
{\bf Proof}\ \ By symmetry, we only need to deal with the last identity. Calculating
modulo $Z(\langle x\rangle^{G_{\beta}})$ we see that 
\begin{eqnarray*}
    [x,x^{abc},x^{c}] & = & [x,x^{abc},[x,c]] \\
            & = & [x,x^{ab}[x^{ab},c],[x,c]] \\
            & \stackrel{\tiny \mbox{L3.3}}{=} & [x,x^{ab},[x,c]] \\
            & = & [x,x^{ab},x^{c}] \\
            &  \stackrel{\tiny \mbox{L3.1}}{=} & [x^{a},x^{b},x^{c}].
\end{eqnarray*}
This finishes the proof. $\Box$
\begin{lemm} $\gamma_{4}(K)\leq Z(\langle x\rangle^{G_{\beta}})$. 
\end{lemm}
{\bf Proof}\ \ Calculations modulo $Z(\langle x\rangle^{G_{\beta}})$ show that 
\begin{eqnarray*}
  1 & \stackrel{\tiny \mbox{L3.3}}{=} & [x,[x^{ab},c],[x,x^{ab}]] \\
                & = & [x,x^{ab}x^{abc},[x,x^{ab}]] \\
                & \stackrel{\tiny \mbox{L3.3}}{=} & [x,x^{abc},[x,x^{ab}]] \\
                & \stackrel{\tiny \mbox{L3.1}}{=} & [x,x^{abc},[x^{a},x^{b}]] \\
                & = & [x,x^{abc},(x^{a}x^{b})^{2}] \\
                & = & [x,x^{abc},x^{a}x^{b}]^{2}[x,x^{abc},x^{a}x^{b},x^{a}x^{b}].
\end{eqnarray*}
By Lemma 3.2 and the presentation for the largest sandwich group generated by $3$ involutions, we know that
$K$ is nilpotent of class at most $4$ and that $\gamma_{3}(K)^{2}=\{1\}$. We also know that $[x^{c},x^{a},x^{b},x^{b}]=[x^{b},x^{c},x^{a},x^{a}]=1$. By Lemma 3.4 we have that $[x,x^{abc},x^{a}x^{b}]^{2}\in \gamma_{3}(K)^{3}=1$ and thus
\begin{eqnarray*}
 1 & = & [x,x^{abc},x^{a}x^{b},x^{a}x^{b}] \\
 1 & \stackrel{\tiny \mbox{L3.4,L3.2}}{=} & [x,x^{abc},x^{a},x^{b} ][x,x^{abc},x^{b},x^{a}][x,x^{abc},x^{a},x^{a}]
[x,x^{abc},x^{b},x^{b}] \\
  & \stackrel{\tiny \mbox{L3.4}}{=} & [x^{b},x^{c},x^{a},x^{b}][x^{c},x^{a},x^{b},x^{a}][x^{b},x^{c},x^{a},x^{a}]
	[x^{c},x^{a},x^{b},x^{b}] \\
    & = & [[x^{b},x^{c}],[x^{a},x^{b}]][[x^{c},x^{a}],[x^{b},x^{a}]].
\end{eqnarray*}
By symmetry we thus have that modulo $Z(\langle x\rangle^{G_{\beta}})$ we have
   $$[[x^{a},x^{b}],[x^{b},x^{c}]]=[[x^{b},x^{c}],[x^{c},x^{a}]]=[[x^{c},x^{a}],[x^{a},x^{b}]].$$
By Lemma 3.2 we know that $K$ is nilpotent of class at most $4$ and thus these three elements all commute with
$x^{a},x^{b},x^{c}$. They of course all commute with $x$ as well. By Lemma 3.1 the first element $[[x^{a},x^{b}],
[x^{b},x^{c}]]$ commutes with $b$. It follows that it commutes with $x^{ab},x^{bc}$. As $[[x^{b},x^{c}],[x^{c},x^{a}]]$ and $[[x^{c},x^{a}],[x^{a},x^{b}]]$ are equal to $[[x^{a},x^{b}],[x^{b},x^{c}]]$ modulo $Z(\langle x\rangle^{G_{\beta}})$, it follows that they also commute with $x^{ab},x^{bc}$. By symmetry the three elements
all commute with $x^{ca}$ as well. Then as $[[x^{a},x^{b}],[x^{b},x^{c}]]$ commutes with $x^{ac}$ and $b$, it commutes
with $x^{abc}$. The same is then true for $[[x^{b},x^{c}],[x^{c},x^{a}]]$ and $[[x^{c},x^{a}],[x^{a},x^{b}]]$. From the presentation for the $3$ generators sandwich group and the fact that $K$ is nilpotent of class at most $4$, we know that $\gamma_{4}(K)=\langle [[x^{a},x^{b}],[x^{b},x^{c}]],\,[[x^{b},x^{c}],[x^{c},x^{a}]],\,[[x^{c},x^{a}],
[x^{a},x^{b}]]\rangle$ and we have thus shown that $\gamma_{4}(K)\leq Z(\langle x\rangle^{G_{\beta}})$. $\Box$
\begin{prop} $\langle x\rangle^{G_{\beta}}$ is nilpotent of class at most $4$.
\end{prop}
{\bf Proof}\ \ It suffices to show that $x\in Z_{4}(\langle x\rangle^{G_{\beta}})$. As $x$ commutes with
$x^{a},x^{b},x^{c}$, it suffices to show that $[x,x^{abc}],[x,x^{ab}],[x,x^{bc}],[x,x^{ca}]$ are in $Z_{3}(\langle x\rangle^{G_{\beta}})$. By symmetry it suffices to show that $[x,x^{abc}],[x,x^{ab}]$ are in $Z_{3}(\langle x\rangle^{G_{\beta}})$. Again by symmetry, it then suffices to show that 
      $$[x,x^{abc},x^{c}], [x,x^{abc},x^{ab}],[x,x^{ab},x^{c}],[x,x^{ab},x^{abc}],[x,x^{ab},x^{ac}]$$
are in $Z_{2}(\langle x\rangle^{G_{\beta}})$. By Lemma 3.4 and Lemma 3.1, we have modulo $Z(\langle x\rangle^{G_{\beta}})$ that 
          $$[x,x^{abc},x^{c}]=[x^{a},x^{b},x^{c}]=[x,x^{ab},x^{c}].$$
By Lemma 3.3 we know that $\gamma_{5}(H(c))\leq Z(\langle x\rangle^{G_{\beta}})$ and thus module $Z(\langle x\rangle^{G_{\beta}})$, we have (using also the fact that $[x^{ab},c]$ commutes with $x^{ab}$)
\begin{eqnarray*}
     [x,x^{abc},x^{ab}] & = & [x,[x^{ab},c],x^{ab}] \\
                         & = & [x,x^{ab},[x^{ab},c]] \\
                         & = & [x,x^{ab},x^{abc}].
\end{eqnarray*}
By Lemma 3.1 we also have
\begin{eqnarray*}
     [x,x^{ab},x^{abc}] & = & [x,x^{ab},x^{c}]^{ab} \\
   \mbox{}  [x,x^{ab},x^{ac}] & = & [x,x^{ab},x^{c}]^{a}.
\end{eqnarray*}
It thus only remains to show that $[x,x^{ab},x^{c}]\in Z_{2}(\langle x\rangle^{G_{\beta}})$. Using the fact 
that $\gamma_{5}(H(c))\leq Z(\langle x\rangle^{G_{\beta}})$ we see that the commutator of $[x,x^{ab},x^{c}]=
[x,x^{ab},[x,c]]$ with $x,x^{ab}$ and $c$ is in $Z(\langle x\rangle^{G_{\beta}})$. Then using 
Lemma 3.5 we know that $\gamma_{4}(K)\leq Z(\langle x\rangle^{G_{\beta}})$ and thus (using Lemma 3.1) we see
that the commutator of $[x,x^{ab},x^{c}]=[x^{a},x^{b},x^{c}]$ with $x^{a},x^{b},x^{c}$ is in $Z(\langle x\rangle^{G_{\beta}})$. Thus we have seen that the commutator of $[x,x^{ab},x^{c}]$ with $x,x^{a},x^{b},x^{ab}$ and $x^{c},
x^{ac},x^{bc},x^{abc}$ is in $Z(\langle x\rangle^{G_{\beta}})$. It follows that $[x,x^{ab},x^{c}]$ is $Z_{2}(\langle x\rangle^{G_{\beta}})$ and this finishes the proof. $\Box$
\begin{theo} $G_{\beta}$ is finite.
\end{theo}
{\bf Proof}\ \ As $G_{\beta}/\langle x\rangle^{G_{\beta}}$ is abelian of order at most $8$, we have that 
$\langle x\rangle^{G_{\beta}}$ is a finitely generated nilpotent torison group and thus finite. From $G/\langle
x\rangle^{G_{\beta}}$ and $\langle x\rangle^{G_{\beta}}$ being finite, it follows that $G_{\beta}$ is finite. 
$\Box$ \\ \\ 
Having determined that the group $G_{\beta}$ is finite, one can obtain the following power commutator presentation
for it. In particular the group has order $2^{28}$. Let 
     $$t(a)=[[x^{c},x^{a}],[x^{a},x^{b}]],\,t(b)=[[x^{a},x^{b}],[x^{b},x^{c}]],\,t(c)=[[x^{b},x^{c}],[x^{c},x^{a}]]$$
and
    $$y(a)=[x,x^{bc},x^{a}],\,y(b)=[x,x^{ca},x^{b}],\,y(c)=[x,x^{ab},x^{c}].$$
\mbox{}\\
\underline{Generators} \\ \\
$X_{1}:\ \ b_{1}=[t(a),b],\ 
b_{2}=t(a),\ b_{3}=t(b)$ \\
$\mbox{\ \ \ \ \ \ \ \ }b_{4}=y(b),\ b_{5}=y(b)^{a},\
b_{6}=y(a),\ b_{7}=y(a)^{b}, b_{8}=y(a)^{c},$ \\ \\
$X_{2}:\ \ b_{9}=[x,x^{ab}][x^{c},x^{abc}],\ b_{10}=[x,x^{ab}],\ 
b_{11}=[x,x^{bc}][x^{a},x^{abc}],$ \\
\mbox{\ \ \ \ \ \ \ \ }$b_{12}=[x,x^{bc}],\ b_{13}=[x,x^{ac}][x^{b},x^{abc}],\ 
b_{14}=[x,x^{ac}],\ b_{15}=[x,x^{abc}][x^{c},x^{ab}]$, \\
\mbox{\ \ \ \ \ \ \ \ }$b_{16}=[x,x^{abc}][x^{a},x^{bc}],\ 
b_{17}=[x,x^{abc}]$, \\ \\
$X_{3}:\ \ b_{18}=x,\ b_{19}=x^{a},\  b_{20}=x^{b},\ b_{21}=x^{c},\ 
b_{22}=x^{ab},\ b_{23}=x^{ca}$, \\
\mbox{\ \ \ \ \ \ \ \ }$b_{24}=x^{bc}, b_{25}=x^{abc}$, \\ \\
$X_{4}:\ \ b_{26}=a,\ b_{27}=b,\ b_{28}=c$. \\ \\
\underline{Relations} \\ \\
$b_{1}^{2}=\ldots =b_{28}^{2}=1$. \\ \\
$b_{2}^{b_{27}}=b_{2}b_{1},\ b_{2}^{b_{28}}=b_{2}b_{1},\ 
b_{3}^{b_{26}}=b_{3}b_{1},\ b_{3}^{b_{28}}=b_{3}b_{1}$ \\
$b_{4}^{b_{19}}=b_{4}b_{2},\ b_{4}^{b_{21}}=b_{4}b_{3}b_{2},\ 
b_{4}^{b_{22}}=b_{4}b_{2}b_{1},
\ b_{4}^{b_{23}}=b_{4}b_{3},\ b_{4}^{b_{24}}=b_{4}b_{3}b_{2}b_{1}$, \\
$b_{4}^{b_{25}}=b_{4}b_{3},\ b_{4}^{b_{26}}=b_{5},\ b_{4}^{b_{28}}=b_{8}b_{6}b_{4}$ \\ 
$b_{5}^{b_{18}}=b_{5}b_{2},\ b_{5}^{b_{20}}=b_{5}b_{2}b_{1},\ b_{5}^{b_{21}}=
b_{5}b_{3}b_{1},\ b_{5}^{b_{23}}=b_{5}b_{3}b_{2}b_{1}$, \\
$b_{5}^{b_{24}}=b_{5}b_{3}b_{1},\ b_{5}^{b_{25}}=b_{5}b_{3}b_{2},\ 
b_{5}^{b_{26}}=b_{4},\ b_{5}^{b_{28}}=b_{8}b_{6}b_{5}$, \\
$b_{6}^{b_{20}}=b_{6}b_{3},\ b_{6}^{b_{21}}=b_{6}b_{3}b_{2},\ b_{6}^{b_{22}}=
b_{6}b_{3}b_{1},\ b_{6}^{b_{23}}=b_{6}b_{3}b_{2}b_{1},\ 
b_{6}^{b_{24}}=b_{6}b_{2}$, \\
$b_{6}^{b_{25}}=b_{6}b_{2},\  
b_{6}^{b_{27}}=b_{7},\ b_{6}^{b_{28}}=b_{8}$, \\
$b_{7}^{b_{18}}=b_{7}b_{3},\ b_{7}^{b_{19}}=b_{7}b_{3}b_{1},\ 
b_{7}^{b_{21}}=b_{7}b_{2}b_{1},\ 
b_{7}^{b_{23}}=b_{7}b_{2}b_{1}$ \\
$b_{7}^{b_{24}}=b_{7}b_{3}b_{2}b_{1},\ b_{7}^{b_{25}}=b_{7}b_{3}b_{2},\ 
b_{7}^{b_{27}}=b_{6},\ b_{7}^{b_{28}}=b_{8}b_{7}b_{6}$, \\
$b_{8}^{b_{18}}=b_{8}b_{3}b_{2},\ b_{8}^{b_{19}}=b_{8}b_{3}b_{2}b_{1},\ 
b_{8}^{b_{20}}=b_{8}b_{2}b_{1},\ b_{8}^{b_{22}}=b_{8}b_{2}b_{1}$, \\
$b_{8}^{b_{24}}=b_{8}b_{3}b_{1},\ b_{8}^{b_{25}}=b_{8}b_{3}, \
b_{8}^{b_{27}}=b_{8}b_{7}b_{6},\ 
b_{8}^{b_{28}}=b_{6}$, \\ \\
$b_{9}^{b_{12}}=b_{9}b_{1},\ 
b_{9}^{b_{14}}=b_{9}b_{1},\ b_{9}^{b_{18}}=b_{9}b_{6}b_{4},$ \\
$b_{9}^{b_{19}}=b_{9}b_{6}b_{5},\
b_{9}^{b_{20}}=b_{9}b_{7}b_{4},\
b_{9}^{b_{21}}=b_{9}b_{6}b_{4}$, \\
$b_{9}^{b_{22}}=b_{9}b_{7}b_{5},\ b_{9}^{b_{23}}=
b_{9}b_{6}b_{5},\ 
b_{9}^{b_{24}}=b_{9}b_{7}b_{4}$, \\
$b_{9}^{b_{25}}=b_{9}b_{7}b_{5}$, \\ 
$b_{10}^{b_{11}}=b_{10}b_{1},\ b_{10}^{b_{12}}=b_{10}b_{3},\ 
b_{10}^{b_{13}}=b_{10}b_{1},\ b_{10}^{b_{14}}=b_{10}b_{2}$, \\
$b_{10}^{b_{16}}=b_{10}b_{1},\ 
b_{10}^{b_{17}}=b_{10}b_{3}b_{2},\ b_{10}^{b_{21}}=b_{10}b_{6}b_{4}$, \\
$b_{10}^{b_{23}}=b_{10}b_{6}b_{5},\ 
b_{10}^{b_{24}}=b_{10}b_{7}b_{4},\ 
b_{10}^{b_{25}}=b_{10}b_{7}b_{5}$, \\
$b_{10}^{b_{28}}=b_{10}b_{9}$, \\
$b_{11}^{b_{14}}=b_{11}b_{1},\ b_{11}^{b_{18}}=b_{11}b_{6},\ 
b_{11}^{b_{19}}=b_{11}b_{6}$, \\
$b_{11}^{b_{20}}=b_{11}b_{7},\ b_{11}^{b_{21}}=b_{11}b_{8},\ 
b_{11}^{b_{22}}=b_{11}b_{7},\ b_{11}^{b_{23}}=b_{11}b_{8},\ 
b_{11}^{b_{24}}=b_{11}b_{8}b_{7}b_{6}$, \\
$b_{11}^{b_{25}}=b_{11}b_{8}b_{7}b_{6}$, \\
$b_{12}^{b_{13}}=b_{12}b_{1},\ b_{12}^{b_{14}}=b_{12}b_{3}b_{2},\ 
b_{12}^{b_{15}}=b_{12}b_{1}$, \\
$b_{12}^{b_{17}}=b_{12}b_{2},\ b_{12}^{b_{19}}=b_{12}b_{6},\ 
b_{12}^{b_{22}}=b_{12}b_{7},\ b_{12}^{b_{23}}=b_{12}b_{8},\ 
b_{12}^{b_{25}}=b_{12}b_{8}b_{7}b_{6}$, \\
$b_{12}^{b_{26}}=b_{12}b_{11}$, \\
$b_{13}^{b_{18}}=b_{13}b_{4},\ b_{13}^{b_{19}}=b_{13}b_{5},\ 
b_{13}^{b_{20}}=b_{13}b_{4}$ \\
$b_{13}^{b_{21}}=b_{13}b_{8}b_{6}b_{4},\ b_{13}^{b_{22}}=b_{13}b_{5},\ 
b_{13}^{b_{23}}=b_{13}b_{8}b_{6}b_{5},\ 
b_{13}^{b_{24}}=b_{13}b_{8}b_{6}b_{4}$ \\ 
$b_{13}^{b_{25}}=b_{13}b_{8}b_{6}b_{5}$, \\
$b_{14}^{b_{15}}=b_{14}b_{1},\ b_{14}^{b_{16}}=b_{14}b_{1},\ 
b_{14}^{b_{17}}=b_{14}b_{3},\ 
b_{14}^{b_{20}}=b_{14}b_{4}$, \\
$b_{14}^{b_{22}}=b_{14}b_{5},\ 
b_{14}^{b_{24}}=b_{14}b_{8}b_{6}b_{4},\ 
b_{14}^{b_{25}}=b_{14}b_{8}b_{6}b_{5},\ 
b_{14}^{b_{27}}=b_{14}b_{13}$, \\
$b_{15}^{b_{18}}=b_{15}b_{6}b_{4}b_{3}b_{2},\
b_{15}^{b_{19}}=b_{15}b_{6}b_{5}b_{3}b_{2}b_{1}$, \\
$b_{15}^{b_{20}}=b_{15}b_{7}b_{4}b_{3}b_{2}b_{1},\ 
b_{15}^{b_{21}}=b_{15}b_{6}b_{4}b_{3}b_{2},\ 
b_{15}^{b_{22}}=b_{15}b_{7}b_{5}b_{3}b_{2},\ 
b_{15}^{b_{23}}=b_{15}b_{6}b_{5}b_{3}b_{2}b_{1}$, \\ 
$b_{15}^{b_{24}}=b_{15}b_{7}b_{4}b_{3}b_{2}b_{1},\ 
b_{15}^{b_{25}}=b_{15}b_{7}b_{5}b_{3}b_{2}$, \\
$b_{16}^{b_{18}}=b_{16}b_{6}b_{2},\ 
b_{16}^{b_{19}}=b_{16}b_{6}b_{2},\ 
b_{16}^{b_{20}}=b_{16}b_{7}b_{2}b_{1}$, \\ 
$b_{16}^{b_{21}}=b_{16}b_{8}b_{2}b_{1},\ 
b_{16}^{b_{22}}=b_{16}b_{7}b_{2}b_{1},\ 
b_{16}^{b_{23}}=b_{16}b_{8}b_{2}b_{1},\ 
b_{16}^{b_{24}}=b_{16}b_{8}b_{7}b_{6}b_{2}$, \\ 
$b_{16}^{b_{25}}=b_{16}b_{8}b_{7}b_{6}b_{2}$, \\
$b_{17}^{b_{19}}=b_{17}b_{6}b_{2},\ 
b_{17}^{b_{20}}=b_{17}b_{4}b_{3},\ 
b_{17}^{b_{21}}=b_{17}b_{6}b_{4}b_{3}b_{2},\ 
b_{17}^{b_{22}}=b_{17}b_{7}b_{5}b_{3}b_{2}$, \\
$b_{17}^{b_{23}}=b_{17}b_{8}b_{6}b_{5}b_{3},\ 
b_{17}^{b_{24}}=b_{17}b_{8}b_{7}b_{6}b_{2},\ 
b_{17}^{b_{26}}=b_{17}b_{16},\ b_{17}^{b_{27}}=b_{17}b_{16}b_{15}$, \\
$b_{17}^{b_{28}}=b_{17}b_{15}$, \\ \\
$b_{18}^{b_{22}}=b_{18}b_{10},\ b_{18}^{b_{23}}=b_{18}b_{14},\ 
b_{18}^{b_{24}}=b_{18}b_{12},\ b_{18}^{b_{25}}=b_{18}b_{17},\ 
b_{18}^{b_{26}}=b_{19}$, \\
$b_{18}^{b_{27}}=b_{20},\ b_{18}^{b_{28}}=b_{21}$, \\
$b_{19}^{b_{20}}=b_{19}b_{10},\ b_{19}^{b_{21}}=b_{19}b_{14},\ 
b_{19}^{b_{24}}=b_{19}b_{17}b_{16},\ b_{19}^{b_{25}}=b_{19}b_{12}b_{11},\
b_{19}^{b_{26}}=b_{18}$, \\ 
$b_{19}^{b_{27}}=b_{22},\ b_{19}^{b_{28}}=b_{23}$, \\
$b_{20}^{b_{21}}=b_{20}b_{12},\ b_{20}^{b_{23}}=b_{20}b_{17}b_{16}b_{15},\ 
b_{20}^{b_{25}}=b_{20}b_{14}b_{13},\ b_{20}^{b_{26}}=b_{22},\ 
b_{20}^{b_{27}}=b_{18}$, \\ 
$b_{20}^{b_{28}}=b_{24}$, \\
$b_{21}^{b_{22}}=b_{21}b_{17}b_{15},\ b_{21}^{b_{25}}=b_{21}b_{10}b_{9},\
b_{21}^{b_{26}}=b_{23},\ b_{21}^{b_{27}}=b_{24},\ b_{21}^{b_{28}}=b_{18}$, \\
$b_{22}^{b_{23}}=b_{22}b_{12}b_{11},\ b_{22}^{b_{24}}=b_{22}b_{14}b_{13},\ 
b_{22}^{b_{26}}=b_{20},\ b_{22}^{b_{27}}=b_{19},\ 
b_{22}^{b_{28}}=b_{25}$, \\
$b_{23}^{b_{24}}=b_{23}b_{10}b_{9},\ b_{23}^{b_{26}}=b_{21},\ 
b_{23}^{b_{27}}=b_{25},\ b_{23}^{b_{28}}=b_{19}$, \\
$b_{24}^{b_{26}}=b_{25},\ b_{24}^{b_{27}}=b_{21},\ 
b_{24}^{b_{28}}=b_{20}$, \\
$b_{25}^{b_{26}}=b_{24},\ b_{25}^{b_{27}}=b_{23},\
b_{25}^{b_{28}}=b_{22}$. \\ 
\subsection{Sandwich groups with commutativity graph $\gamma$}
Let $G_{\gamma}=\langle a,b,x,y\rangle$ be a sandwich group generated by involutions whose commutativity
graph is \\ \\ \\
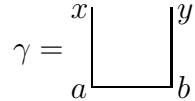
\begin{picture}(100,50)(-120,-50)
\put(45,-25){\line(1,0){30}}
%\put(30,20){\line(3,5){30}}
%\put(90,20){\line(-3,5){30}}
\put(77,-29){$b$}
\put(37,-29){$a$}
\put(15,-14){$\gamma=$}
\put(77,1){$y$}
\put(37,1){$x$}
%\put(145,5){\line(1,-1){30}}
%\put(145,-25){\line(1,1){30}}
\put(45,5){\line(0,-1){30}}
\put(75,5){\line(0,-1){30}}
\end{picture}
\normalsize
\mbox{}
\begin{prop}
$G_{\gamma}$ is finite. 
\end{prop}
{\bf Proof}\ \ We have that $\langle y,b,y^{a}\rangle$ is abelian and thus $H=\langle
y,b,y^{a},x\rangle$ is a homomorphic image of $G_{\beta}$. As $H\unlhd G_{\gamma}$ and
$G/H=\langle aH\rangle$, it follows that $G_{\gamma}$ is finite. $\Box$ \\ \\
Again one can then obtain a power-conjugation presentation of $G_{\gamma}$. This turns out to be the
following one. In particular $G_{\gamma}$ has order $2^{20}$. \\ \\
\underline{Generators} \\ \\
$X_{1}:\ \ e_{1}=[x,b,[y,a],x,y,[x,b]],\ 
e_{2}=[x,b,[y,a],y,x,[y,a]]$, \\
$\mbox{\ \ \ \ \ \ \ \ }e_{3}=[x,b,[y,a],x,y,x],\ 
e_{4}=[x,b,[y,a],y,x,y]$, \\ 
$\mbox{\ \ \ \ \ \ \ \ }e_{5}=[x,b,[y,a],x,y],\ 
e_{6}=[x,b,[y,a],y,x]$, \\
$\mbox{\ \ \ \ \ \ \ \ }e_{7}=[x,b,[y,a],x,\ 
e_{8}=[x,b,[y,a],y]$, \\ 
$\mbox{\ \ \ \ \ \ \ \ }e_{9}=[x,b,y,x],\ e_{10}=
[x,[y,a],y]$ \\ \\
$X_{2}:\ \ e_{11}=[x,b,[y,a]],\ e_{12}=[x,b,y],\ 
e_{13}=[x,[y,a]]$, \\ 
\mbox{\ \ \ \ \ \ \ \ }$e_{14}=[x,y]$ \\ \\
$X_{3}:\ \ e_{15}=[x,b],\ e_{16}=[y,a],\ 
e_{17}=x,\ e_{18}=y$ \\ \\
$X_{4}:\ \ e_{19}=a,\ e_{20}=b$. \\ \\
\underline{Relations} \\ \\
$e_{1}^{2}=e_{2}^{2}=\ldots =e_{20}^{2}=1$. \\ \\
$e_{3}^{e_{20}}=e_{3}e_{1},\ e_{4}^{e_{19}}=e_{4}e_{2},\ 
e_{5}^{e_{15}}=e_{5}e_{1},\ e_{5}^{e_{17}}=e_{5}e_{3},$ \\
$e_{6}^{e_{16}}=e_{6}e_{2},\ e_{6}^{e_{18}}=e_{6}e_{4},\ 
e_{7}^{e_{12}}=e_{7}e_{1},\ e_{7}^{e_{14}}=e_{7}e_{3},\ 
e_{7}^{e_{18}}=e_{7}e_{5}$,\\
$e_{8}^{e_{13}}=e_{8}e_{2},\ e_{8}^{e_{14}}=e_{8}e_{4},\ e_{8}^{e_{17}}=e_{8}e_{6}$ \\
$e_{9}^{e_{11}}=e_{9}e_{1},\ e_{9}^{e_{13}}=e_{9}e_{3},\ 
e_{9}^{e_{16}}=e_{9}e_{5}e_{4}e_{2},\ e_{9}^{e_{19}}=
e_{9}e_{7}e_{6}e_{1}$, \\
$e_{10}^{e_{11}}=e_{10}e_{2},\ e_{10}^{e_{12}}=e_{10}e_{4},\
e_{10}^{e_{15}}=e_{10}e_{6}e_{3}e_{1},\ e_{10}^{e_{20}}=
e_{10}e_{8}e_{5}e_{2}$, \\
$e_{11}^{e_{14}}=e_{11}e_{6}e_{5}e_{4}e_{3},\ e_{11}^{e_{17}}=e_{11}e_{7},\ 
e_{11}^{e_{18}}=e_{11}e_{8}$, \\
$e_{12}^{e_{13}}=e_{12}e_{6}e_{5}e_{4}e_{3},\ 
e_{12}^{e_{16}}=e_{12}e_{8},\ e_{12}^{e_{17}}=e_{12}e_{9},\ 
e_{12}^{e_{19}}=e_{12}e_{11}e_{8}$, \\
$e_{13}^{e_{15}}=e_{13}e_{7},\ e_{13}^{e_{18}}=e_{13}e_{10},\ 
e_{13}^{e_{20}}=e_{13}e_{11}e_{7}$, \\
$e_{14}^{e_{15}}=e_{14}e_{9},\ e_{14}^{e_{16}}=e_{14}e_{10},\ 
e_{14}^{e_{19}}=e_{14}e_{13}e_{10},\ e_{14}^{e_{20}}=e_{14}e_{12}e_{9}$, \\
$e_{15}^{e_{16}}=e_{15}e_{11},\ e_{15}^{e_{18}}=e_{15}e_{12},\ 
e_{16}^{e_{17}}=e_{16}e_{13}$, \\
$e_{17}^{e_{18}}=e_{17}e_{14},\ e_{17}^{e_{20}}=e_{17}e_{15},\ 
e_{18}^{e_{19}}=e_{18}e_{16}$. 
\mbox{} \\ \\ \\
\section{Proof of the main result}
In this final section we give a proof of Theorem 1.4. As we saw in the introduction, this reduces to ${\mathcal Q}(2,2)$ and
${\mathcal Q}(5,1)$. For the proof of ${\mathcal Q}(2,2)$ we will need the preliminary work from Section 2 and 3 whereas the proof of ${\mathcal Q}(5,1)$ will
be modeled on the approach in [11]. 

\subsection{Proof of ${\mathcal Q}(2,2)$}
Let $G$ be a group with a left $3$-Engel element $x$ of order $2$. Suppose $\langle x\rangle^{G}$ has no elements of order $8$. The aim is to show that $x\in HP(G)$. That is we want to show that $E=\langle x\rangle^{G}$ is locally nilpotent. For this we only need to show that $x\in HP(E)$. Without loss of generality we can thus replace $G$ by $E$ and assume that $G$ has no elements of order $8$. As groups of exponent $4$ are locally finite [10], it suffices to show
that ${\mathcal R}_{G}(2,2)$ holds. \\ \\
Let $g\in G$ be an element of order $4$. Let $H=\langle a,c,b,d\rangle=\langle x,x^{g},x^{g^{2}},x^{g^{3}}\rangle$. By Proposition 1.3 it suffices to show that $H$ is nilpotent and this will be our aim.
\begin{lemm}
If $h\in G$ is an involution and $y$ is a conjugate of $x$, then $[y,y^{h}]=1$. 
\end{lemm}
{\bf Proof}\ \ We have $\langle y,h\rangle =\langle y,y^{h}\rangle \rtimes \langle h\rangle$ where
$\langle y,y^{h}\rangle$ is nilpotent of class at most $2$. Hence $\langle y,h\rangle$ is a $2$-group. As
$G$ has no element of order $8$, we have that 
     $$1=(hy)^{4}=h^{4}(y^{h^{3}}y^{h^{2}}y^{h}y)=(yy^{h})^{2}=[y^{h},y].$$
This finishes the proof. $\Box$
\begin{lemm} If $y,z$ are two conjugates of $x$ that commute and $h$ is an involution in $G$, then 
   $$[y^{h},z]=[y,z^{h}].$$
\end{lemm}
{\bf Proof}\ \  By Lemma 4.1 we have that $y,z^{h}$ commute with $y^{h},z$ and 
 $$\langle h,y,z\rangle=\langle y,z,y^{h},z^{h}\rangle \rtimes \langle h\rangle=(\langle y,z^{h}\rangle 
\cdot \langle y^{h},z\rangle) \rtimes \langle h\rangle.$$
As $\langle y,z^{h}\rangle$ is nilpotent of class at most $2$ it is clear that this group is finite of order at most $8^{2}\cdot 2=128$. As there is no element of order $8$, we then must have
$$
   1=(h(yz))^{4}=[(yz)^{h}(yz)]^{2}=(y^{h}z)^{2}(z^{h}y)^{2}=[y^{h},z][z^{h},y].$$
Hence the result. $\Box$ \\ \\
{\bf Remark}\ \ By Lemma 4.1 it follows in particular that $[x,x^{g^{2}}]=[x^{g},x^{g^{3}}]=1$. Thus if
we pick any three of the generators $x,x^{g},x^{g^{2}},x^{g^{3}}$ then we know that two of them must commute. \\
\begin{lemm}
Let $u,v,t$ be three of the generators of $x,x^{g},x^{g^{2}},x^{g^{3}}$ where $u$ and $v$ commute. Then $\langle
u,v,t\rangle$ is nilpotent of class at most $3$.
\end{lemm}
{\bf Proof}\ \ We have $[u,t,v,v]=[v^{[u,t]},v]=1$ by Lemma 4.1 as $1=[u,t^{2}]=[u,t]^{2}$. Similarly we have
that $[v,t,u,u]=[u^{[t,v]},u]=1$. From the presentation of the largest $3$-generator sandwich group, generated by
involutions, we thus know that $\langle u,v,t\rangle$ is nilpotent of class at most $4$. Finally
$[\textbf{\textbf{}}[t,u],[t,v]]=[t^{u},t^{v}]=[t^{u},(t^{u})^{uv}]=1$ by Lemma 4.1 as $(uv)^{2}=1$. Again by the presentation
for the largest $3$-generator sandwich group we see that $\langle u,v,t\rangle$ is nilpotent of class at most
$3$. $\Box$
\begin{prop} Let $T=\langle a,b,x,y\rangle$ be a subgroup of $G$ that is a sandwich group of type 
\mbox{}\\ \\ \\
\begin{picture}(100,50)(-120,-50)
\put(45,-25){\line(1,0){30}}
%\put(30,20){\line(3,5){30}}
%\put(90,20){\line(-3,5){30}}
\put(77,-29){$b$}
\put(37,-29){$a$}
\put(15,-14){$\gamma=$}
\put(77,1){$y$}
\put(37,1){$x$}
%\put(145,5){\line(1,-1){30}}
%\put(145,-25){\line(1,1){30}}
\put(45,5){\line(0,-1){30}}
\put(75,5){\line(0,-1){30}}
\end{picture}
\mbox{}\\
where $a,b,x,y$ are involutions. Then $T$ is nilpotent of class at most $3$ and $T^{4}=\{1\}$. 
\end{prop}
{\bf Proof}\ \ We already know from Section 3 that the group $T$ is a finite $2$-group and as there are no elements
of order $8$ we must have $T^{4}=\{1\}$. In the following calculations we use the presentation for the largest
group of type $\gamma$ generated by involutions. With a slight abuse of notation we will use $e_{1},e_{2},\ldots
,e_{20}$ for the values of this free group in $T$ under the natural homomorphism. As $(xa)^{2}=(yb)^{2}=1$ we can deduce from Lemma 4.1 that 
    $$1=[y^{x},y^{a}]=[y,x,[y,a]]=[e_{14},e_{16}]=e_{10}$$
and 
    $$1=[x^{y},x^{b}]=[x,y,[x,b]]=[e_{14},e_{15}]=e_{9}.$$
It is easy to see that from this it follows that $e_{1}=\ldots =e_{8}=1$. Also
$$1=(byxa)^{4}=[(xa)^{by}xa]^{2}=(x^{by}a^{y}xa)^{2}=(x^{by}aa^{y}x)^{2}.$$
Notice that $1=[x^{y},x^{b}]^{y}=[x,x^{by}]$. Now we thus get
  $$1=[(xx^{by})(aa^{y})]^{2}=[(aa^{y})^{xx^{by}}(aa^{y})=[aa^{y},xx^{by}]=[a,xx^{by}]^{a^{y}}[a^{y},xx^{by}]
	=[a,x^{by}][a^{y},x].$$
Hence 
   $$1=[a^{y},x^{b}][a,x^{y}]=[a^{y},x]^{b}[a,[x,y]]=[e_{16},e_{17}]^{e_{20}}[e_{19},e_{14}]=e_{11}.$$
From the presentation one now observes that $e_{1}=\ldots =e_{11}=1$ implies that $T$ is nilpotent of class
at most $3$. $\Box$ \\ \\
We now return to the main task of this section. Let $H=\langle a,b,c,d\rangle=\langle x,x^{g^{2}},x^{g},x^{g^{3}}\rangle$. We know that $[a,b]=[c,d]=1$. We want to show this group is a finite $2$-group. 
\begin{lemm} We have that $[[a,c],[b,d]],\,[[a,d],[b,c]]\leq Z_{2}(H)$.
\end{lemm}
{\bf Proof}\ \ We first show that $\langle [a,c],\,[b,d]\rangle$ and $\langle [a,d],\,[b,c]\rangle$ are nilpotent
of class at most $2$. We first turn to the first one. Calculating in $\langle a,x,b,y\rangle =\langle a,a^{c},
b,b^{d}\rangle$ we get a group of type $\gamma$. Hence using the presentation for the largest such group and the fact from the last proposition that this group has class at most $3$, we see that
   $$[a,c,[b,d],[a,c]]=[ax,by,ax]=[a,y,x][x,y,a]=[e_{16},e_{17}][e_{14},e_{19}]=1.$$
By symmetry we also have that $[[a,c],[b,d]]$ commutes with $[b,d]$. Similarly one sees that $\langle [a,d],\,[b,c]\rangle$ is nilpotent of class at most $2$. Notice that $1=[a,c]^{2}=[b,d]^{2}=[a,d]^{2}=[b,c]^{2}$ and thus $[a,c]=[c,a],[b,d]=[d,b],[a,d]=[d,a],[b,c]=[c,b]$. From what we have seen above we also know that
$$ [[a,c],[b,d]]=[[b,d],[a,c]],\ [[a,d],[b,c]]=[[b,c],[a,d]].$$
This symmetry implies that in order to show that $[[a,c],[b,d]],[[a,d],[b,c]]$ are in $Z_{2}(H)$, it suffices that
their commutator with $a$ is in $Z(H)$. Calculating again in $\langle a,x,b,y\rangle = \langle a,a^{c},b,b^{d}\rangle
$, we see that 
     $$[a,c,[b,d],a]=[ax,by,a]=[x,y,a]=[a^{c},b^{d},a].$$
As $[c,d]=[a,b]=1$, we have by Lemma 4.2 that $[a^{c},b^{d}]=[a^{d},b^{c}]$. Thus 
\begin{equation}
  [a,c,[b,d],a]=[a^{c},b^{d},a]=[a^{d},b^{c},a]=[a,d,[b,c],a].
\end{equation}
Notice the joint element commutes with $a$ and $b$ as $\langle a,a^{c},b,b^{d}\rangle$ is nilpotent
of class at most $3$. Next notice that in $\langle a,x,b,y\rangle=\langle a,a^{c},b,b^{d}\rangle$ we have
     $$[a,c,[b,d],a]=[x,y,a]=[e_{14},e_{19}]=e_{13}=[e_{17},e_{16}]=[x,[y,a]]=[ax,[by,a]]=[a,c,[b,d,a]].$$
Notice also that $[b,d,a,a]=[by,a,a]=[y,a,a]=1$ and thus $1=[b,d,a^{2}]=[b,d,a]^{2}[b,d,a,a]=[b,d,a]^{2}$. By Lemmma
4.1 we thus have that 
      $$\langle c,c^{a},c^{[b,d,a]},c^{a[b,d,a]}\rangle$$
is abelian. Hence 
\begin{eqnarray*}
  [a,c,[b,d,a]]^{c} & = & [a,c,[b,d,a][b,d,a,c]] \\
                    & = & [a,c,[b,d,a,c]]\cdot [a,c,[b,d,a]]^{[b,d,a,c]} \\
    \mbox{}                & = & [a,c,[b,d,a]].
\end{eqnarray*}
Thus $[a,c,[b,d],a]=[a,c,[b,d,a]]$ commutes with $c$ and by symmetry $[a,d,[b,c],a]$ commutes with $d$. From (4) we thus now know that $[a,c,[b,d],a]=[a,d,[b,c],a]$ commutes with $a,b,c,d$ and is thus in $Z(G)$. $\Box$
\begin{prop}
We have that $\langle a,b,c,d\rangle$ is finite.
\end{prop}
{\bf Proof}\ \ By Lemma 4.1 we know that $\langle c,c^{a},c^{b},c^{ab}\rangle$ is abelian and thus
$[c,a,b]$ commutes with $c$. Clearly $[c,a,b]$ commutes with $a,b$. By symmetry $[c,b,a]$ commutes with $a,b,c$ and
thus $\langle a,b,c\rangle$ is nilpotent of class at most $3$. From this an the fact that $[a,b]=1$, it follows that
$[c,a,b]=[c,b,a]$. In order to show that $[c,a,b]\in Z_{3}(\langle a,b,c,d\rangle)$ it then only remains to see that
$[c,a,b,d]=[c,b,a,d]\in Z_{2}(\langle a,b,c,d\rangle)$. It clearly commutes with $c$ as $[c,d]=1$. As $[c,a,b]^{2}=1$
it follows from Lemma 4.1 that it commutes with $d$. Then
    $$[c,a,b,d,b]=[c,a,b,d]^{-1}[c,a,b,d[d,b]]=[c,a,b,d]^{-1}[c,a,b,[d,b]][c,a,b,d]^{[d,b]}.$$
By the last lemma we know that $[c,a,b]^{[d,b]}=[[c,a][[c,a],[d,b]],b]=[c,a,b]$ modulo $Z(\langle a,b,c,d\rangle)$. 
Thus $[c,a,b,d,b]=[c,a,b,d]^{-1}[c,a,b,d]^{[d,b]}=[c,a,b,d]^{-1}[c,a,b,d]=1$ modulo $Z(\langle a,b,c,d\rangle)$.
By symmetry $[c,a,b,d,a]=[c,b,a,d,a]\in Z(\langle a,b,c,d\rangle)$. We have thus shown that $[c,a,b]\in Z_{3}(\langle
a,b,c,d\rangle)$. By symmetry this is true for any commutator of weight $3$ in $a,b,c,d$. Hence $\langle a,b,c,d
\rangle$ is nilpotent of class at most $6$. $\Box$

\subsection{Proof of ${\mathcal Q}(5,1)$}
Let $G$ be a group with a left $3$-Engel element $x$ of order $5$ and suppose furthermore that $H=\langle x\rangle^{G}$
has no element of order $25$. The aim is to show that $x$ is in the locally nilpotent radical of $G$. The proof
of this will be modeled on [11], where this is proved under the stronger hypothesis that $G$ is of exponent $5$. A key ingredient is a certain variant of a similar result from [12]. Before stating it we recall some terminology from [5]. We say that a group $\langle a,b,c\rangle$ is of type $(r,s,t)$ if $\langle a,b\rangle$, $\langle a,c\rangle$ and
$\langle b,c\rangle$ are of class at most $r,s$ and $t$ respectively. 
\\
\begin{prop}
Let $K=\langle a,b,c\rangle$ be a subgroup of $H$ that is of type $(1,2,3)$ and where $a,b,c$ are of order $5$. Suppose
furthermore that $c$ is a left $3$-Engel element of $H$ and that $[b,c,c]=1$. The $\langle a,b,c\rangle$ is 
nilpotent of class at most $4$ and of exponent $5$.
\end{prop}
{\bf Proof}\ \ From the proof of Proposition 2.2 in [11] (see Step 2), we know that $\langle c\rangle^{K}$ is nilpotent of class at most $3$. As $c$ is of order $5$ it follows in particular that $\langle c\rangle^{K}$ is a $5$-group. Now
$K/\langle c\rangle^{K}$ is abelian of order dividing $25$. Thus 
$\langle c\rangle^{K}$ is a finitely generated nilpotent $5$-group and thus a finite $5$-group. Hence $K$ is a finite $5$-group. As $K\leq H$, it has no elements of order $25$ and thus must be of exponent $5$. The fact that $K$ is nilpotent of class
at most $4$ now follows from Proposition 2.3 in [11]. $\Box$ \\ \\
In order to show that $x$ is in the locally nilpotent radical 
of $G$ is suffices to show that $\langle x\rangle^{G}$ is 
locally nilpotent. It thus suffices to prove the following result
that is again a variant of the corresponding result in [11] for
groups of exponent $5$. That proof was also modelled on a similar result in [12].
\begin{prop} Let $k$ be a positive integer and let $a_{1},\ldots ,a_{k}$ be conjugates of $x$. Then $A=\langle a_{1},\ldots 
,a_{k}\rangle$ is nilpotent of class at most $k$ and of exponent $5$. Furthermore $\langle a_{i}\rangle^{A}$ is abelian
for $i=1,\ldots ,k$. 
\end{prop}
{\bf Proof}\ \ From the proof of Theorem 3.1 in [11], we know that this result holds when 
$A$ is of exponent $5$. Thus it suffices to show that $A$ is nilpotent as then $A$ is a $5$-group and the assumption that $H$ has no elements of order $25$ implies then that $A$ is of exponent $5$. \\ \\
 We have that the case $k=2$ holds by the asumption that $x$ is a left $3$-Engel element and the case $k=3$ follows from the fact that $3$-generator sandwich groups are nilpotent. Now suppose that $k\geq 3$. Let 
$u=[a_{1},a_{2},\ldots ,a_{k-2}]$.  Then the subgroup $\langle a_{k-1},a_{k-1}^{u},a_{k}\rangle$ is generated by $3$-conjugates of
$x$ and is thus nilpotent of class at most $3$. By this and the fact that any two conjugates generate a subgroup of class at most $2$, it follows that 
  $$[a_{1},a_{2},\ldots ,a_{k-1},a_{k},a_{k}]=[a_{k-1}^{-u}a_{k-1},a_{k},a_{k}]=1$$
and 
  $$[a_{k},_{3}[a_{1},a_{2},\ldots ,a_{k-1}]]=[a_{k},_{3}a_{k-1}^{-u}a_{k-1}]=1.$$
We thus have the following identities which hold for any conjugates $a_{1},a_{2},\ldots ,a_{k}$ of $x$ and for any $k\geq 3$. 
\begin{eqnarray}
    [a_{1},a_{2},\ldots, a_{k-1},a_{k},a_{k}] & = & 1, \\ \nonumber
 [a_{k},[a_{1},a_{2},\cdots ,a_{k-1}],
   [a_{1},a_{2},\ldots ,a_{k-1}],
	[a_{1},a_{2},\ldots ,a_{k-1}]] & = & 1.
\end{eqnarray}
We now proceed with the induction step. Let $k\geq 4$ and suppose that the result is true for all smaller values of $k$. We first show that if $1\leq r\leq k$, then 
\begin{equation}
 [[a_{1},a_{2},\ldots ,a_{r}],[a_{1},a_{k},a_{k-1},\ldots ,a_{r+1}]]=[a_{1},a_{2},\ldots ,a_{k},a_{1}]^{(-1)^{k-r}}.
\end{equation}
This is obvious when $r=k$. Now consider the
case $r=k-1$. Let $u=[a_{1},\ldots ,a_{k-1}]$. By the induction hypothesis and (4) we have that $\langle a_{1},u,a_{k}\rangle$ is
of type $(1,2,3)$ and satifies the condition for Proposition 4.7. Hence $\langle a_{1},u,a_{k}\rangle$ is nilpotent of class at most $4$. Using the fact that $u$ commutes with $a_{1}$ and the first identity in (4) one sees easily that all commutators of weight $(2,1,1)$ and $(1,1,2)$ in $a_{1},u,a_{k}$ are trivial. The only commutators that one needs to consider are $[u,a_{k},a_{1},a_{1}]$ and $[u,a_{k},a_{1},a_{k}]$ but as $[u,[a_{k},a_{1},a_{1}]]=
[u,[a_{1},a_{k},a_{k}]]=1$ we get by expanding these that 
\begin{eqnarray*}
    1 & = & [u,a_{k},a_{1},a_{1}], \\
    1 & = & [u,a_{k},a_{1},a_{k}].
\end{eqnarray*}
From this one sees that $[[u,[a_{1},a_{k}]]=[u,a_{k},a_{1}]^{-1}$ that gives us identity (5) when $r=k-1$. This argument also tells us that 
  $$[[a_{1},a_{k},\ldots ,a_{3}],[a_{1},a_{2}]]=[a_{1},a_{k},\ldots ,a_{2},a_{1}]^{-1}$$
and thus
   $$[a_{1},[a_{1},a_{k},\ldots ,a_{2}]]=
	[[a_{1},a_{2}],[a_{1},a_{k},\ldots ,a_{3}]]^{-1}$$
that shows that the case $r=1$ follows if it holds for $r=2$. To establish (5) it is thus sufficient to show that for $2\leq r\leq k-2$ we have 
 $$[[a_{1},a_{2},\ldots ,a_{r}],[a_{1},a_{k},\ldots ,a_{r+1}]]=
[[a_{1},a_{2},\ldots ,a_{r+1}],[a_{1},a_{k},\ldots ,a_{r+2}]]^{-1}.$$
Let $u=[a_{1},a_{2},\ldots ,a_{r}]$ and
$v=[a_{1},a_{k},\ldots ,a_{r+2}]$. By the induction hypothesis we have that $u$ and $v$ commute and that $\langle u,a_{r+1}\rangle,
\langle v,a_{r+1}\rangle$ are nilpotent of
class at most $2$. Thus $\langle u,v,a_{r+1}\rangle$ is of type $(1,2,2)$. From induction hypothesis we also know that $u^{5}=v^{5}=1$. By Section 2.1.2. in [11] we know that the group $\langle u,v,a_{r+1}\rangle$ is nilpotent and thus of exponent $5$. The presentation given for this group in Section 2.1.2 in [11] shows that the group is then nilpotent of class at most $3$. Thus $[u,[v,a_{r+1}]]=[u,a_{r+1},v]^{-1}$ as required. This establishes (5). \\ \\
We want to show that $A$ is nilpotent. We will show that $A$ is nilpotent of class at most $k+1$. The rest will then follow from the fact that $A$ is of exponent $5$ and [11]. \\ \\
Consider a commutator $c=[b_{1},b_{2},\ldots ,b_{k+1}]$ where $b_{1},\ldots ,b_{k+1}\in
\{a_{1},\ldots ,a_{k}\}$. We want to show that $c\in Z(A)$. By induction $c=1$ unless 
$\{b_{1},\ldots ,b_{k}\}=\{a_{1},\ldots ,a_{k}\}$. Also by (4) we have that $c=1$ if $b_{k}=b_{k+1}$. So there is no loss of generality in assuming that $b_{k+1}=a_{1},
b_{k}=a_{k}$ and that $\{b_{1},\ldots ,b_{k-1}\}=\{a_{1},\ldots ,a_{k-1}\}$. Then, using the inductive hypothesis, we see that $[b_{1},b_{2},\ldots ,b_{k-1}]$ can be expressed as a product $u_{1}u_{2}\cdots u_{r}$ where each $u_{i}$ is a commutator of the form 
$[a_{1},a_{\sigma(2)},a_{\sigma(3)},\ldots ,a_{\sigma(k-1)}]^{\epsilon_{i}}$ for some permutation $\sigma$ of $\{2,3,\ldots ,k-1\}$ and where $\epsilon_{i}=\pm 1$. So 
  $$c=[b_{1},\ldots ,b_{k+1}]=
	[u_{1}\cdots u_{r},a_{k},a_{1}]=
	[\prod_{i=1}^{r}[u_{i},a_{k}]^{u_{i+1}u_{i+2}\cdots u_{r}},a_{1}].$$
Now the inductive hypothesis implies that $u_{1},u_{2},\ldots ,u_{r}$ commute with $a_{1}$. So $c$ is the product of conjugates of the commutators $[u_{1},a_{k},a_{1}],\ldots ,[u_{r},a_{k},a_{1}]$. Notice also that $[u_{i},a_{k}]$ commutes with $a_{1}$ if and only if $[u_{i}^{-1},a_{k}]$ commutes with $a_{1}$. To show that $c\in Z(A)$ it thus clearly suffices to show that $[a_{1},a_{2},\ldots ,a_{k},a_{1}]\in Z(A)$. \\ \\
So consider $d=[a_{1},a_{2},\ldots ,a_{k},a_{1},a_{i}]$, where $1\leq i\leq k$. If $i=1$ then $d=1$ by (4). If $i=k$, let $u=[a_{1},a_{2},\ldots ,a_{k-1}]$. Then, using the induction hypothesis, $u$ is of order $5$, $\langle a_{1},u,a_{k}\rangle$ is of type $(1,2,3)$ and satisfies the conditions given in Proposition 4.7. It is thus nilpotent of class at most $4$ (and then exponent $5$). Thus
   $$1=[u,[a_{1},a_{k},a_{k}]]=[u,a_{k},a_{1},a_{k}]^{-2}.$$
This implies that $[u,a_{k},a_{1},a_{k}]=1$ and thus $d=1$ when $i=k$. Now let $1<i<k$. To show that $d=1$, it suffices by (5) to show that $[u,a_{i},v,a_{i}]=1$ when $u=[a_{1},a_{2},\ldots ,a_{i-1}]$ and $v=[a_{1},a_{k},a_{k-1},\ldots ,a_{i+1}]$. Now by the induction hypothesis $\langle u,v,a_{i}\rangle$ is of type $(1,2,3)$ with $u^{5}=v^{5}=1$ satisfies the criteria from Proposition 4.7. Thus it is nilpotent of class at most $4$. Hence 
again 
     $$1=[u,[v,a_{i},a_{i}]]=[u,a_{i},v,a_{i}]^{-2}$$
that implies that $[a_{1},a_{2},\ldots ,a_{k},a_{1}]$ commutes with $a_{i}$. This finishes the proof that $A$ is nilpotent and thus the inductive step. $\Box$

\end{document}